\documentclass[12pt]{amsart}

\usepackage{epsf}
\usepackage{amssymb}

\newcommand{\Gr}{\mbox{\it Grass}}

\newcommand{\llra}{\relbar\joinrel\longrightarrow}

\newcommand{\cpp}{\stackrel{c_p}{\llra}}

\newtheorem{lem}{Lemma}[section]
\newtheorem{thm}[lem]{Theorem}

\newtheorem{cor}[lem]{Corollary}

\newtheorem{alg}[lem]{Algorithm}

\begin{document}

\title[Schubert Polynomials and Orders on ${\mathcal S}_n$]{Identities of
Structure Constants for Schubert Polynomials
and Orders on ${\mathcal S}_n$}

\author{Nantel Bergeron \and Frank Sottile}

\address{Department of Mathematics and Statistics\\
        York University\\
        North York, Ontario M3J 1P3\\
	CANADA}
\email[Nantel Bergeron]{bergeron@@mathstat.yorku.ca}
\address{Department of Mathematics\\
        University of Toronto\\
        100 St.~George Street\\
	Toronto, Ontario  M5S 3G3\\
	CANADA}
\email[Frank Sottile]{sottile@math.toronto.edu}
\date{1 May 1997}
\thanks{First author supported in part by an NSERC grant}
\thanks{Second author supported in part by NSERC grant  OGP0170279 \ 
and NSF grant DMS-9022140}
\thanks{To appear in FPSAC'97 (Formal Power Series and Algebraic
        Combinatorics, Vienna, 1997) conference abstracts}

\maketitle

\def\baselinestretch{.90}
\begin{center}
\begin{minipage}[t]{5.8in}\footnotesize
{\sc R\'esum\'e.}
L'analogue des coefficients de Littlewood-Richardson pour les polyn\^omes
de Schubert est reli\'e \`a l'\'enum\'eration de chaines dans
l'ordre partiel de Bruhat de ${\mathcal S}_n$. Pour mieux comprendre ce
lien, nous le raffinons en r\'eduisant le probl\`eme \`a des sous-ordres
partiels li\'es aux sous-groupes paraboliques du groupe sym\'etrique. Nous
montrons ici certaines identit\'es g\'eom\'etriques reliant ces coefficients
entre eux et, pour la plupart de ces identit\'es, nous montrons des
r\'esultats combinatoires compagnons pour les chaines dans l'ordre de Bruhat.
Nous esp\'erons que la compr\'ehension du lien entre les chaines et les
coefficient permettra la d\'eduction des identit\'es g\'eom\'etriques \`a
partir
des identit\'es combinatoires.
Dans ces travaux, nous donnons: un nouvel ordre partiel gradu\'e sur le
groupe sym\'etrique, des r\'esultats sur l'\'enum\'eration de chaines dans
l'ordre de Bruhat, ainsi qu'une formule pour une grande vari\'et\'e de
sp\'ecialisations des variables pour les polyn\^omes de Schubert.
\bigskip

{\sc Summary.}
For Schubert polynomials, the analogues of Littlewood-Richard\-son
coefficients are expected to be related to the enumeration of chains in
the Bruhat order on ${\mathcal S}_n$.
We refine this expectation in terms of certain suborders on the
symmetric group associated to parabolic subgroups.
Our main results are a number of new identities among these
coefficients. 
For many of these identities, there is a companion result about the
Bruhat order which we expect would imply the identity, were it known
how to express these coefficients in terms of the Bruhat order.
Our analysis leads to a new graded partial order on the symmetric group,
results on the enumeration of chains in the Bruhat order, the determination
of many of these constants, and formulas for a large class of
specializations of the variables in a Schubert polynomial.
\end{minipage}
\end{center}

\def\baselinestretch{1.00}

\section*{Introduction}
Extending work of Demazure~\cite{Demazure} and of Bernstein, Gelfand, and
Gelfand~\cite{BGG},  in 1982 Lascoux and
Sch\"utzenberger~\cite{Lascoux_Schutzenberger_polynomes_schubert}
defined remarkable polynomial representatives for Schubert classes in
the cohomology of a flag manifold, which they called 
Schubert polynomials.
For each permutation $w$ in ${\mathcal S}_\infty$, the group of permutations
of ${\mathbb N} :=\{1,2,\ldots\}$ which fix all but finitely many numbers,
there is a Schubert polynomial ${\mathfrak S}_w\in{\mathbb Z}[x_1,x_2,\ldots]$.
The collection of all Schubert polynomials forms an additive basis for
this polynomial ring.
Thus the identity
	\begin{equation}\label{eq:structure}
	{\mathfrak S}_u \cdot {\mathfrak S}_v \quad 
	=\quad \sum_w c^w_{u\, v} {\mathfrak S}_w
	\end{equation}
defines integral {\em structure constants} $c^w_{u\, v}$ for the ring of
polynomials with respect to its Schubert basis.
Littlewood-Richardson coefficients are
a special case of the $c^w_{u\, v}$ as 
any Schur symmetric polynomial  is a Schubert polynomial.
The  $c^w_{u\, v}$  are non-negative integers:
Evaluating a Schubert polynomial at certain Chern classes
gives a Schubert class in the cohomology of the flag manifold.
Hence, $c^w_{u\, v}$ enumerates the flags in a
suitable triple intersection of Schubert varieties.
It is an open problem to give a combinatorial interpretation or a bijective
formula for these constants.

All known formulas express $c^w_{u\, v}$ in terms of chains in
the Bruhat order.
For instance, the Littlewood-Richardson
rule~\cite{Littlewood_Richardson}, a special case, may be expressed in
this form.  
Other formulas for these constants, particularly Monk's
formula~\cite{Monk}, Pieri-type
formulas~\cite{Lascoux_Schutzenberger_polynomes_schubert,%
sottile_pieri_schubert,Winkel_multiplication,Kirillov_Maeno}, and the other
formulas of~\cite{sottile_pieri_schubert}, are all of this form. 
For quantum Schubert
polynomials~\cite{Fomin_Gelfand_Postnikov,Ciocan_partial},
the Pieri-type formulas~\cite{Ciocan_partial,Postnikov}
are also of this form.

We present a number of new identities among the $c^w_{u\,v}$
which are consistent with the expectation that they can be expressed in
terms of chains in the Bruhat order.
In addition, we give a formula (Theorem~\ref{thm:skew_shape}) for
many of the  $c^w_{u\,v}$ when ${\mathfrak S}_v$ is a symmetric polynomial.
These identities impose stringent conditions on 
any combinatorial interpretation for the coefficients.
They also point to some potentially beautiful combinatorics
once such an interpretation is known.

These  results are expanded on and proven in a
manuscript, ``{S}chubert polynomials, the {B}ruhat order, and the
geometry of flag manifolds''~\cite{bergeron_sottile_symmetry}.
For a background on Schubert polynomials, we recommend the original
papers~\cite{Lascoux_Schutzenberger_polynomes_schubert,%
Lascoux_Schutzenberger_structure_de_Hopf,%
Lascoux_Schutzenberger_symmetry,%
Lascoux_Schutzenberger_interpolation,%
Lascoux_Schutzenberger_schub_LR_rule,%
Lascoux_Schutzenberger_operators},
the survey~\cite{Lascoux_historique}, or the
book~\cite{Macdonald_schubert}.
For their relation to geometry, we recommend the
book~\cite{Fulton_tableaux}.

\section{Chains in the Bruhat order}\label{sec:suborders}

Let $(a,b)$ denote the transposition interchanging $a<b$.
The {\em Bruhat order} $\leq $ on the symmetric group ${\mathcal S}_n$ 
is defined by its covers:
$$
u\ \lessdot\ u(a,b) \quad \Longleftrightarrow\quad
\ell(u) + 1 \ =\ \ell(u(a,b)).
$$
It also appears as the index of summation in Monk's
formula~\cite{Monk}: 
$$		
{\mathfrak S}_u \cdot {\mathfrak S}_{(k,\,k{+}1)} \quad =\quad 	
{\mathfrak S}_u \cdot (x_1+\cdots+x_k)
 \quad =\quad \sum  {\mathfrak S}_{u(a, b)},
$$	
the sum over all $a\leq k<b$ where $\ell(u(a,b))=\ell(u)+1$.

This suggests the following notion, which appeared
in~\cite{Lascoux_Schutzenberger_symmetry}. 
A {\em coloured chain} is a (saturated) chain in the
Bruhat order together with 
an element of $\{a,a+1,\ldots,b-1\}$ for each cover 
$u\lessdot u(a,b)$ in that chain.
Let $I$ be any subset of ${\mathbb N}$.
An {\em $I$-chain} is a coloured chain whose colours are chosen from the
set $I$. 
If $u\leq w$ in the Bruhat order, let $f^w_u(I)$ count the 
$I$-chains from $u$ to $w$.
Iterating Monk's formula, we obtain:
$$
	\left(\sum_{i\in I}\:{\mathfrak S}_{(i,i{+}1)}\right)^m
	\quad =\quad 
	\sum_v\, f^v_e(I)\; {\mathfrak S}_v.
$$
Multiplying this expression by ${\mathfrak S}_u$, expanding the products
using (\ref{eq:structure}) and Monk's formula, and equating
coefficients of ${\mathfrak S}_w$, we obtain:

\begin{thm}\label{thm:chains}
Let $u,w\in{\mathcal S}_\infty$ and $I\subset{\mathbb N}$.
Then
$$
f^w_u(I)\quad =\quad 
\sum_v\, c^w_{u\,v}\, f^v_e(I).
$$
\end{thm}

This number, $f^v_e(I)$, is non-zero precisely when $v$ is minimal in
its coset $vW_{\overline{I}}$, where $W_{\overline{I}}$ is the parabolic
subgroup~\cite{Bourbaki_Groupes_IV} of ${\mathcal S}_\infty$ generated by the
transpositions $\{(i,i{+}1)\:|\: i\not\in I\}$.
We say that $u$ is comparable to $w$ in the {\em $I$-Bruhat order}
if  there is an $I$-chain from $u$ to $w$.
In \S\S\ref{sec:k-bruhat} and~\ref{sec:Schur_identities}, we consider
this order when $I=\{k\}$.

Any eventual combinatorial interpretation of the constants
$c^w_{u\,v}$ should give a bijective proof of 
Theorem~\ref{thm:chains}.
We expect there will be a combinatorial interpretation of the
following form:
Let $u,v,w\in{\mathcal S}_\infty$, and $I\subset{\mathbb N}$ be such
that $v$ is minimal in $vW_{\overline{I}}$.
(There always is such an $I$.)
Then,  for any $I$-chain $\gamma$ from $e$
to $v$,
$$
	c^w_{u\, v}\quad =\quad\# \left\{
	\mbox{\begin{minipage}[c]{2.15in}
	$I$-chains from
	$u$ to $w$   satisfying 
	\,{\em some} condition \,imposed by \,$\gamma$ \end{minipage}}
	\right\}.
$$

\section{The $k$-Bruhat order}\label{sec:k-bruhat}

The $k$-Bruhat  order, $\leq_k$, is the $\{k\}$-Bruhat order of
\S\ref{sec:suborders}.  
It has another description:

\begin{thm}\label{thm:k-length}
Let $u,w\in {\mathcal S}_\infty$ and $k\in{\mathbb N}$. 
Then $u\leq _k w$ if and only if
\begin{enumerate}
\item[I.] $a\leq k < b$ implies $u(a)\leq w(a)$ and $u(b)\geq w(b)$.
\item[II.] If $\/a<b$, $u(a)<u(b)$, and $w(a)>w(b)$, then 
$a\leq k< b$.
\end{enumerate}
\end{thm}

Considering covers shows conditions I and II are necessary.
Sufficiency follows from a greedy algorithm:

\begin{alg}[Produces a chain in the $k$-Bruhat order]\label{alg:chain}
\mbox{ }

\noindent{\tt input: }Permutations $u,w\in {\mathcal S}_\infty$ 
satisfying conditions I and II of Theorem~\ref{thm:k-length}.

\noindent{\tt output: }A chain in the $k$-Bruhat order from $w$ to $u$.

Output $w$.
While $u\neq w$, do
\begin{enumerate}
\item[1] Choose $a\leq k$ with $u(a)$ minimal subject to $u(a)< w(a)$.
\item[2] Choose $b>k$ with $u(b)$ maximal subject to $w(b)<w(a)\leq u(b)$.
\item[3] (Then $w(a,b)\lessdot_k w$.) Set $w:=w(a,b)$, output $w$.   
\end{enumerate}

At every iteration of\/ {\rm 1}, $u,w$ satisfy conditions I and II of
Theorem~\ref{thm:k-length}. 
Moreover, this algorithm terminates in $\ell(w)-\ell(u)$ iterations and the
sequence of permutations produced is a  chain in the $k$-Bruhat
order from $w$ to $u$.
\end{alg}

Observe that Algorithm~\ref{alg:chain} may be stated in
terms of the permutation $\zeta:=w u^{-1}$:
\medskip

\noindent{\tt input: }{\it A permutation $\zeta\in {\mathcal S}_\infty$.}

\noindent{\tt output: }{\it Permutations
$\zeta,\zeta_1,\ldots,\zeta_m=e$ such that  if $u\leq_k \zeta u$, then
$$
u\ \lessdot_k\ \zeta_{m-1}u\ \lessdot_k\ \cdots\ 
\lessdot_k\ \zeta_1u\ \lessdot_k\ \zeta
u
$$
is a saturated chain in the $k$-Bruhat order.

Output $\zeta$.
While $\zeta\neq e$, do
\begin{enumerate}
\item[1] Choose $\alpha$ minimal subject to $\alpha< \zeta(\alpha)$.
\item[2] Choose $\beta$ maximal subject to 
$\zeta(\beta)<\zeta(\alpha)\leq\beta$. 
\item[3] $\zeta:=\zeta(\alpha,\beta)$,  output $\zeta$.
\end{enumerate}}
To see this is equivalent to Algorithm~\ref{alg:chain}, 
set $\alpha=u(a)$ and $\beta=u(b)$ so that $w(a)=\zeta(\alpha)$ and
$w(b)=\zeta(\beta)$.
Thus 
$w(a, b) = \zeta u(a, b) = \zeta(\alpha,\beta) u$.
This observation is generalised considerably in 
Theorem~\ref{thm:B} ({\em i}) below.

\section{Identities when ${\mathfrak S}_v$ is a Schur
polynomial}\label{sec:Schur_identities} 

The Schur symmetric polynomial $S_\lambda(x_1,\ldots,x_k)$ is the
Schubert polynomial ${\mathfrak S}_{v(\lambda,k)}$, where $v(\lambda,k)$ is
a {\em Grassmannian permutation}, a permutation with unique descent at $k$.
Here $\lambda_{k+1-j} = v(j)-j$ and $v(\lambda,k)$ is minimal in its 
coset $v(\lambda,k) W_{\overline{\{k\}}}$.
Consider the constants $c^w_{u\:v(\lambda,k)}$ which are defined by
the identity:
$$
{\mathfrak S}_u \cdot S_\lambda(x_1,\ldots,x_k) \quad =\quad
\sum_{w} c^w_{u\:v(\lambda,k)}{\mathfrak S}_w.
$$
By Theorem~\ref{thm:chains}, the $c^w_{u\:v(\lambda,k)}$ are related to the
enumeration of chains in the $k$-Bruhat order.
These $c^w_{u\:v(\lambda,k)}$  share many properties with 
Littlewood-Richardson coefficients:

If $\lambda, \mu$, and $\nu$ are partitions with at most $k$ parts, then 
the Littlewood-Richardson coefficients 
$c^\nu_{\mu\,\lambda}$ are defined by the identity
$$
S_\mu(x_1,\ldots,x_k)\cdot S_\lambda(x_1,\ldots,x_k) \quad = 
\quad \sum_\nu c^\nu_{\mu\,\lambda} S_\nu(x_1,\ldots,x_k).
$$
The $ c^\nu_{\mu\,\lambda}$ depend only upon $\lambda$ and
the skew partition $\nu/\mu$.
That is, if $\kappa$ and $\rho$ are partitions with at most $l$ parts
and $\kappa/\rho = \nu/\mu$, then 
for any partition $\lambda$,
$$
c^\nu_{\mu\,\lambda} \quad=\quad
c^\kappa_{\rho\,\lambda}.
$$
Moreover, $c^\kappa_{\rho\,\lambda}$ is the coefficient of 
$S_\kappa(x_1,\ldots,x_l)$ when 
$S_\rho(x_1,\ldots,x_l)\cdot S_\lambda(x_1,\ldots,x_l)$
is expressed as a sum of Schur polynomials.
The order type of the interval in Young's lattice
from $\mu$ to $\nu$ is determined by $\nu/\mu$.

If $u\leq_k w$, let $[u,w]_k$ be the interval between $u$ and $w$ in the
$k$-Bruhat order.
Permutations $\zeta$ and $\eta$ are {\em shape equivalent} if there
exist sets of integers $P=\{p_1<\cdots<p_n\}$ and $Q=\{q_1<\cdots<q_n\}$,
where $\zeta$ (respectively $\eta$) acts as the identity on 
${\mathbb N}- P$ (respectively ${\mathbb N}- Q$), and 
$$
\zeta(p_i)\ =\ p_j \quad \Longleftrightarrow\quad
\eta(q_i)\ =\ q_j.
$$

\begin{thm}[Skew Coefficients]\label{thm:B}
Suppose $u\leq_kw$ and $x\leq_l z$ where
$wu^{-1}$ is shape equivalent to $zx^{-1}$.
Then
\begin{enumerate}
\item[({\em i})] $[u,w]_k\simeq[x,z]_l$.
When $wu^{-1}=zx^{-1}$, this isomorphism is given by $v\mapsto vu^{-1}x$.
\item[({\em ii})] For any partition $\lambda$ with length at most the
minimum of $l$ and $k$, 
$$
c^w_{u\,v(\lambda,k)} = c^z_{x\,v(\lambda,l)}.
$$
\end{enumerate}
\end{thm}

By Theorem~\ref{thm:B} ({\em ii}), we may define the 
{\em skew coefficient} 
$c^\zeta_\lambda$ for $\zeta\in {\mathcal S}_\infty$ and 
$\lambda$ a partition by 
$c^\zeta_\lambda := c^{\zeta u}_{u\,v(\lambda,k)}$ for any 
$u \in {\mathcal S}_\infty$ with $u\leq_k \zeta u$.
(There always is a $u$ and $k$ with $u\leq_k \zeta u$.)
Moreover, $c^\zeta_\lambda$ depends only upon $\lambda$ and the shape
equivalence class of $\zeta$.

Theorem~\ref{thm:B} ({\em i}) is proven
using combinatorial arguments.
The key lemma is that if 
$u\lessdot_k(\alpha,\beta)u\leq_k w$ and $y\leq_k z$ with 
$wu^{-1}=zy^{-1}$, then 
$y\lessdot_k(\alpha,\beta)y\leq_k z$.

The identity of Theorem~\ref{thm:B} ({\em ii}) is proven 
using geometric arguments ({\em cf.}~\S5):
It follows from an equality of homology classes,
which we show by explicitly computing the intersection of two Schubert
varieties in a flag manifold and the image of that intersection under a
projection to a Grassmannian. 

By Theorem~\ref{thm:B} ({\em i}), we may define a partial
order  $\preceq$ on ${\mathcal S}_\infty$ as follows:
Set $\eta\preceq \zeta$ if there exists $u\in {\mathcal S}_\infty$ and $k$
such that $u\leq_k\eta u\leq_k\zeta u$.
This partial order has a rank function defined by 
$|\zeta|:=\ell(\zeta u)-\ell(u)$, whenever $u\leq_k \zeta u$.  
Both the definition of $\preceq$ and Theorem~\ref{thm:B} ({\em i}) 
are illustrated by the following example:

Let $\zeta = (24)(153)$ and $\eta=(35)(174)$.
Then $\zeta$ nd $\eta$  are shape equivalent.
Also 
$21345 \leq_2 45123 = \zeta\cdot 21345$
and
$3215764\leq_3 5273461= \eta\cdot 3215764$.
We illustrate the intervals		%Figure~\ref{fig:new_order} 
$[21342,\,45123]_2$, $[3215764,\,5273461]_3$, and 
$[e,\zeta]_\preceq$ below

$$\epsfxsize=5.in \epsfbox{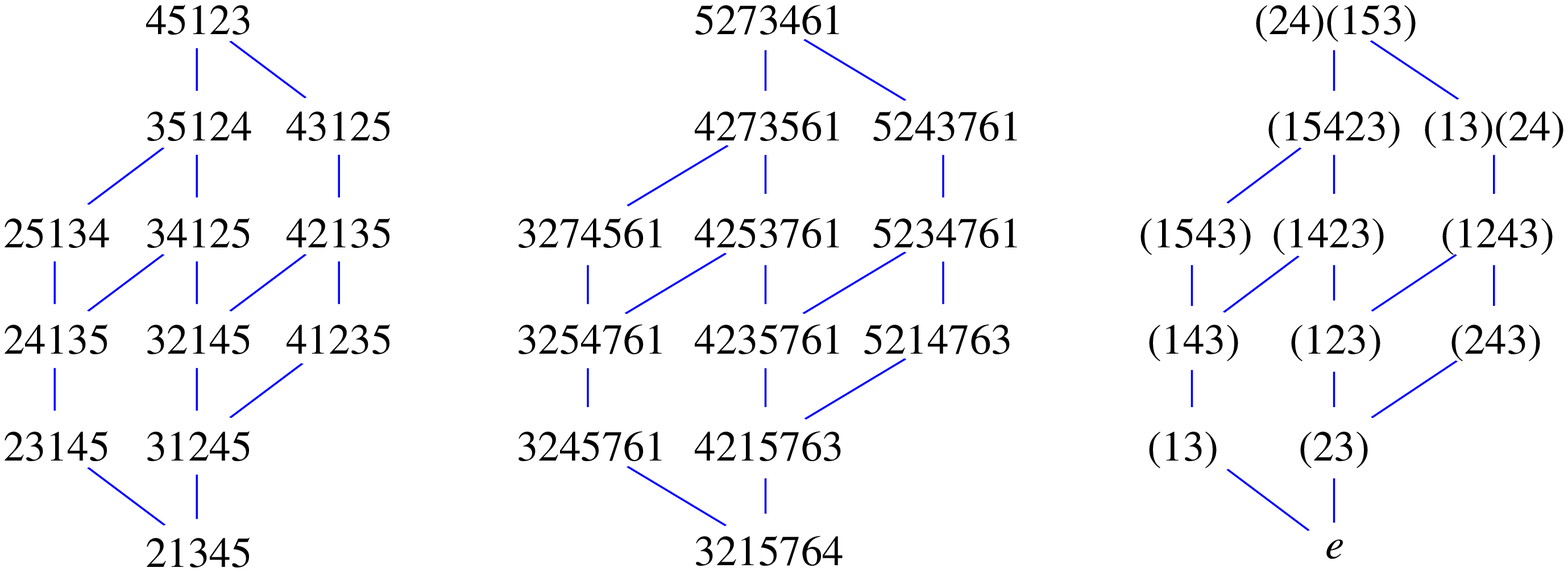}$$

The order and rank function may be defined independent of 
$u$ and $k$:
Define $\mathrm{up}_\zeta:= \{\alpha\;|\; \zeta(\alpha)>\alpha\}$ and 
$\mathrm{down}_\zeta:=\{\alpha\;|\; \zeta(\alpha)<\alpha\}$.
Then $\eta\preceq\zeta$ if and only if 
\begin{enumerate}
\item $\alpha<\eta(\alpha)\Rightarrow \eta(\alpha)\leq \zeta(\alpha)$,
\item $\alpha>\eta(\alpha)\Rightarrow \eta(\alpha)\geq \zeta(\alpha)$,
and
\item If $\alpha,\beta\in\mathrm{up}_\zeta$ or 
$\alpha,\beta\in\mathrm{down}_\zeta$  with $\alpha<\beta$ and
$\zeta(\alpha)<\zeta(\beta)$, then $\eta(\alpha)<\eta(\beta)$.
\end{enumerate}
Similarly, $|\zeta|$ equals the difference of 
$\#\{(\alpha,\beta)\in 
\zeta(\mbox{up}_\zeta)\times\zeta(\mbox{down}_\zeta)\,|\,
\alpha>\beta\}$ and 
$$
\begin{array}{c}
\#\{(a,b)\in \mbox{up}_\zeta\times\mbox{down}_\zeta\,|\,
a>b\} +\ 
 \#\{(a,b)\in \mbox{up}_\zeta\times\mbox{up}_\zeta\,|\,
a>b \mbox{ and  } \zeta(a)<\zeta(b)\}
\\ +\  
\#\{(a,b)\in \mbox{down}\zeta\times\mbox{down}\zeta\,|\,
a>b \mbox{ and } \zeta(a)<\zeta(b)\}.
\end{array}
$$

This new order is preserved by many group-theoretic operations:
For $\zeta\in {\mathcal S}_n$, let $\overline{\zeta}:=w_0\zeta w_0$,
conjugation by the longest element of ${\mathcal S}_n$.
For $P: p_1<p_2<\cdots \subset{\mathbb N}$ and $\zeta\in{\mathcal S}_\infty$,
define the homomorphism 
$\phi_P:{\mathcal S}_\infty \rightarrow {\mathcal S}_\infty$, 
by requiring that 
$\phi_P(\zeta)\in {\mathcal S}_\infty$ act as the identity on ${\mathbb N}-P$ 
and $\phi_P(\zeta)(p_i)=p_{\zeta(i)}$.
Note that $\zeta$ and $\phi_P(\zeta)$ are shape equivalent.

\begin{thm}\label{thm:new_order}
Suppose $\zeta,\eta,\xi\in{\mathcal S}_\infty$.
\begin{enumerate}
\item[({\em i})] 
The restriction of $\preceq$ to Grassmannian permutations $v(\lambda,k)$
gives Young's lattice of partitions with at most $k$ parts.
\item[({\em ii})] 
For $\eta\preceq\zeta$, the map 
$\xi\mapsto \xi\eta^{-1}$ induces
an isomorphism $[\eta,\zeta]_\preceq \stackrel{\sim}{\longrightarrow}
[e,\zeta\eta^{-1}]_\preceq$.
\item[({\em iii})] 
For any infinite set  $P\subset {\mathbb N}$,
$\phi_P:{\mathcal S}_\infty \rightarrow {\mathcal S}_\infty$ is 
an injection of graded posets.
\item[({\em iv})] 
The map $\eta\mapsto \eta\zeta^{-1}$ is an order reversing 
bijection
between $[e,\zeta]_\preceq$ and $[e,\zeta^{-1}]_\preceq$.
\item[({\em v})] 
The homomorphism $\zeta\mapsto \overline{\zeta}$ on ${\mathcal S}_n$ 
induces an
automorphism of $({\mathcal S}_n,\preceq)$.
\end{enumerate}
\end{thm}

\begin{figure}[htb]\label{fig:new_order.S_4}
	$$\epsfxsize=4in \epsfbox{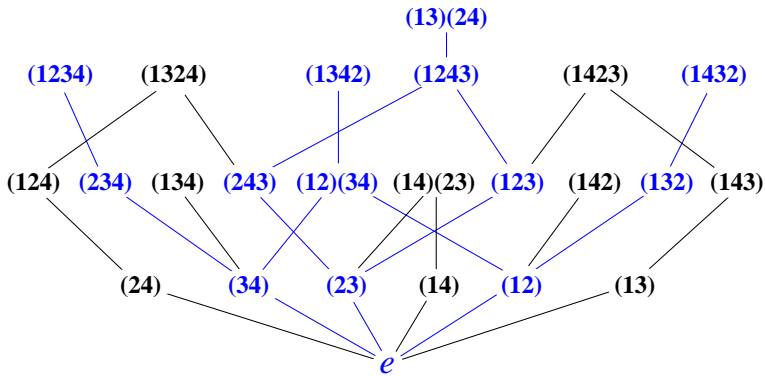}$$
	\caption{ $\preceq$ on  ${\mathcal S}_4$}
	\end{figure}

These properties are easy consequences of the definitions.
This order is studied further in~\cite{bergeron_sottile_order}. 
Figure~\ref{fig:new_order.S_4} shows $\preceq$ on ${\mathcal S}_4$.

Some of these structure constants 
$c^w_{u\,v(\lambda,k)}$ may be expressed in
terms of chains in the Bruhat order.
If $u\lessdot_k u(a,b)$ is a cover in the $k$-Bruhat order, label that edge
of the Hasse diagram with the integer $u(b)$.
The {\em word} of a chain in the $k$-Bruhat order is the sequence of
labels of edges in the chain.

\begin{thm}\label{thm:skew_shape}
Suppose $u\leq_k w$ and $wu^{-1}$  is shape equivalent to
$v(\mu,l)\cdot v(\nu,l)^{-1}$, for some $l$ and partitions
$\mu\subset\nu$. 
Then, for any partition $\lambda$ and standard Young tableau $T$ of
shape $\lambda$, 
$$
c^w_{u\,v(\lambda,k)}\ =\ 
\#\left\{\begin{array}{cc}\mbox{Chains in $k$-Bruhat order from $u$ to
$w$ whose word} \\\mbox{  has recording tableau $T$ under Schensted
insertion} \end{array}\right\}.
$$
\end{thm}

Since a chain in the $k$-Bruhat order is a $\{k\}$-chain in the sense of
\S\ref{sec:suborders}, 
Theorem~\ref{thm:skew_shape} gives a combinatorial proof of
Theorem~\ref{thm:chains} when
$wu^{-1}$ is shape equivalent to a skew partition.
The key step is when $w$ and $u$ are
Grassmannian permutations.
In that case, symmetry of the Schensted algorithm reduces the
theorem to showing the `diagonal word' of a tableau is Knuth-equivalent
to its reading word.
By the {\em diagonal word}, we mean the entries of a tableau read 
first by their diagonal, and then in increasing order in each diagonal.
For instance, this tableau has diagonal word 
$7\,58\,379\,148\,26\,26\,5\,8$.
$$
\epsfxsize=.9in \epsfbox{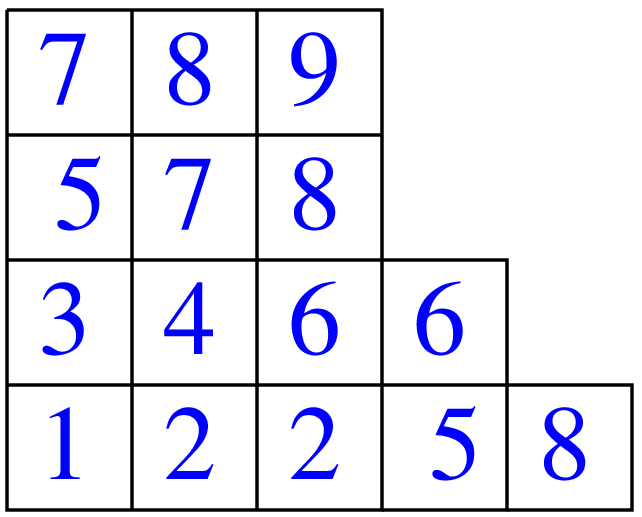}
$$

There are many permutations $u,w$ which do not satisfy the
hypotheses of Theorem~\ref{thm:skew_shape}, but for which the conclusion
holds.
For example, any $u,w$ for which $wu^{-1}=(143652)$ satisfy the
conclusion, but $(143652)$ is not shape equivalent to any skew
partition.
However, some hypotheses are necessary.
Let $u=312645$ and $w=561234$.
Here is the labeled Hasse diagram of $[u,w]_2$:
$$
\epsfxsize=1.9in \epsfbox{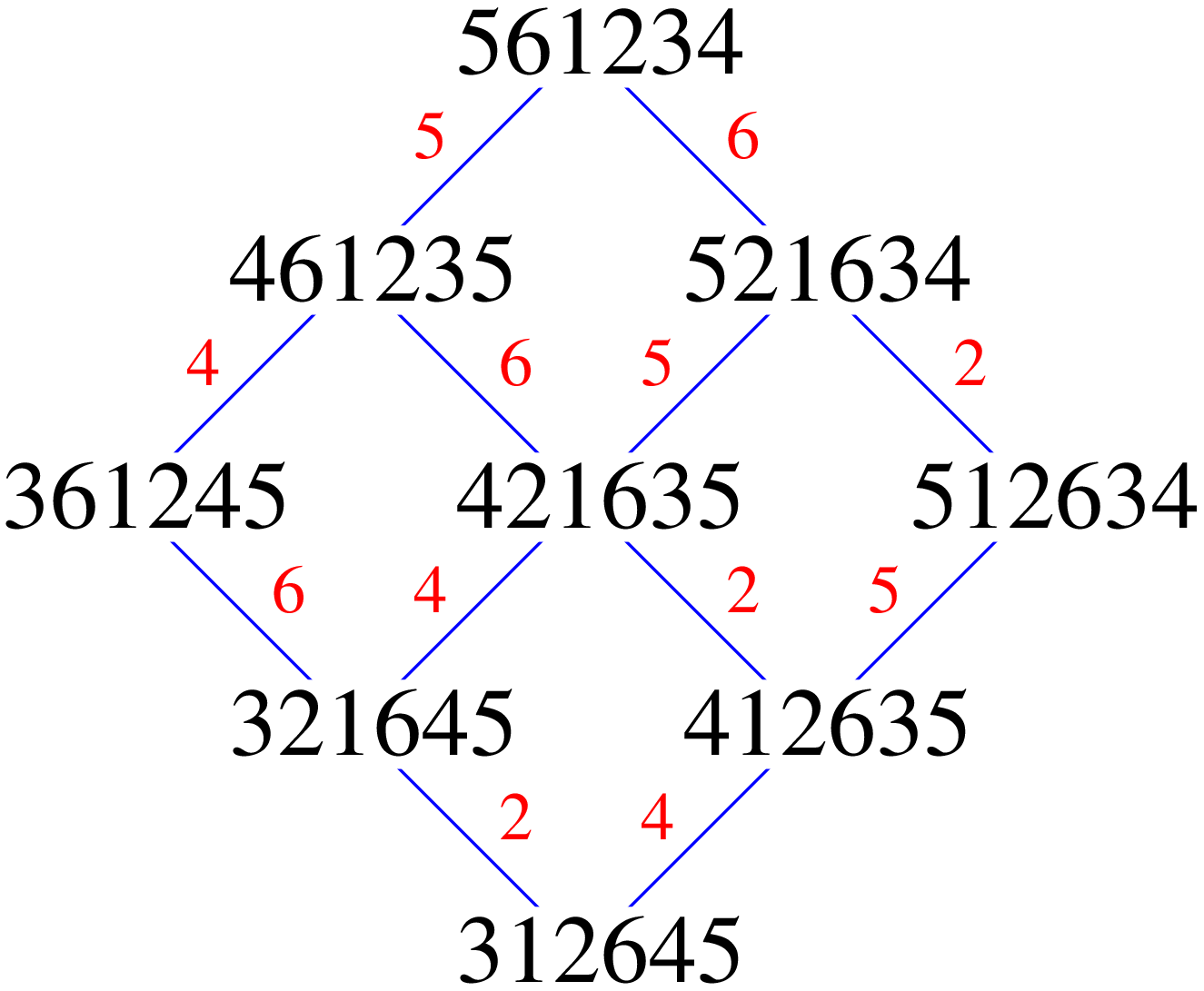}
$$
There are six chains in this interval
from which we obtain these recording tableaux:
$$
\epsfxsize=4.3in \epsfbox{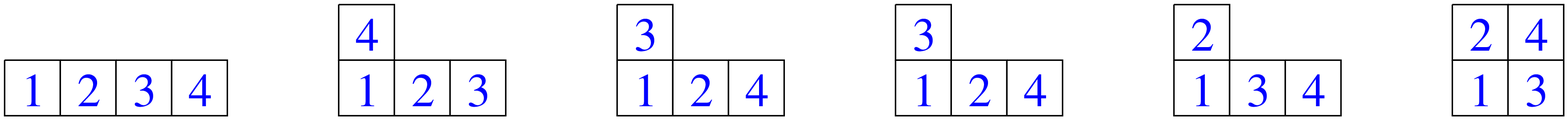}\ \raisebox{6pt}{.}
$$
This list omits the tableau
\begin{picture}(17,17)(0,2)   \put(0,0){\line(0,1){17}} 
\put(0,0){\line(1,0){17}}\put(8.5,0){\line(0,1){17}}
\put(0,8.5){\line(1,0){17}}\put(17,0){\line(0,1){17}}
\put(0,17){\line(1,0){17}}
\put(1.5,1.5){\scriptsize 1} \put(10,1.5){\scriptsize 2}
\put(1.5,10){\scriptsize 3} \put(10,10){\scriptsize 4}
\end{picture},
and the third and fourth tableaux are identical.
\smallskip

We calculate $c^w_{u,v(\lambda,2)}$ using Theorem~5
of~\cite{sottile_pieri_schubert}, or Theorem 4 (1)
of~\cite{sottile_geometry_schubert}: 
$$
c^w_{u\; v(\begin{picture}(12,3) \put(0,0){\line(0,1){3}}
\put(0,0){\line(1,0){12}}      \put(3,0){\line(0,1){3}}
\put(0,3){\line(1,0){12}}      \put(6,0){\line(0,1){3}}
\put(9,0){\line(0,1){3}}      \put(12,0){\line(0,1){3}}
\end{picture}\,,2)}\ =\ 
c^w_{u\; v(\,\begin{picture}(9,6) \put(0,0){\line(0,1){6}}
\put(0,0){\line(1,0){9}}      \put(3,0){\line(0,1){6}}
\put(0,3){\line(1,0){9}}      \put(6,0){\line(0,1){3}}
\put(0,6){\line(1,0){3}}      \put(9,0){\line(0,1){3}} 
\end{picture}\,,2)}\ =\ 
c^w_{u\; v(\,\begin{picture}(6,6)(0,0) 
\put(0,0){\line(0,1){6}}      \put(0,0){\line(1,0){6}}
\put(3,0){\line(0,1){6}}      \put(0,3){\line(1,0){6}}
\put(6,0){\line(0,1){6}}      \put(0,6){\line(1,0){6}}
\end{picture}\,,2)}\ =\ 1.
$$

If a skew Young diagram $\kappa$ is the disjoint union of two incomparable
diagrams $\rho$ and $\theta$, then 
$$
c^\kappa_\lambda\ =\ 
\sum_{\mu,\nu} c^\lambda_{\mu\,\nu}c^\rho_\mu c^\theta_\nu.
$$
Similarly, we say that a product $\zeta\cdot\eta$ is {\em disjoint} if
the two permutations $\zeta$ and $\eta$ have disjoint supports and 
$|\zeta\cdot\eta| = |\zeta|+ |\eta|$.
We have:

\begin{thm}[Disjointness] Suppose $\zeta\cdot\eta$ is disjoint. 
Then
\begin{enumerate}
\item The map $(\alpha,\beta)\mapsto \alpha\cdot\beta$ induces an
isomorphism $[e,\zeta]_\preceq\times[e,\eta]_\preceq 
\stackrel{\sim}{\longrightarrow}[e,\zeta\cdot\eta]_\preceq$.
\item For all $\lambda$, \ 
${\displaystyle c^{\zeta\cdot\eta}_\lambda\ =\ 
\sum_{\mu,\nu} c^\lambda_{\mu\,\nu}c^\zeta_\mu c^\eta_\nu}$.
\end{enumerate}
\end{thm}

The next identity has no analogy with the
classical Littlewood-Richardson coefficients.
Let $(1\,2\,\ldots\, n)$ be the permutation which cyclicly
permutes $[n]$.

\begin{thm}[Cyclic Shift]\label{thm:D}
Suppose $\zeta\in S_n$ and $\eta = \zeta^{(1\,2\,\ldots\,n)}$.
Then, for any partition $\lambda$, 
$c^\zeta_\lambda = c^\eta_\lambda$.
\end{thm}

This is proven using geometric arguments similar to those which
establish Theorem~\ref{thm:B} ({\em ii}).
Combined with Theorem~\ref{thm:chains},
we obtain:

\begin{cor}\label{cor:equal_chains}
If $u\leq_k w$ and $x\leq_k z$ with $w u^{-1},z x^{-1}\in{\mathcal S}_n$
and 
$(wu^{-1})^{(1\,2\,\ldots\,n)} = zx^{-1}$, then 
the two intervals $[u,w]_k$ and $[x,z]_k$ each have the same number of
chains.
\end{cor}

It would be interesting to give a bijective proof of
Corollary~\ref{cor:equal_chains}. 
The two intervals
$[u,w]_k$ and $[x,z]_k$ of Corollary~\ref{cor:equal_chains} are
typically non-isomorphic: 
In ${\mathcal S}_4$, let $u=1234$, $x=2134$, and $v=1324$.
If $\zeta=(1243)$, $\eta=(1423)= \zeta^{(1234)}$, and 
$\xi=(1342)=\eta^{(1234)}$, then 
$$
u \ \leq_2\ \zeta u,\quad x\ \leq_2\  \eta x,\quad\mbox{and}
\quad v\ \leq_2 \ \xi v.
$$
We illustrate the intervals $[u,\zeta u]_2$,  $[x,\eta x]_2$,
and  $[v,\xi v]_2$.
$$\epsfxsize=4in \epsfbox{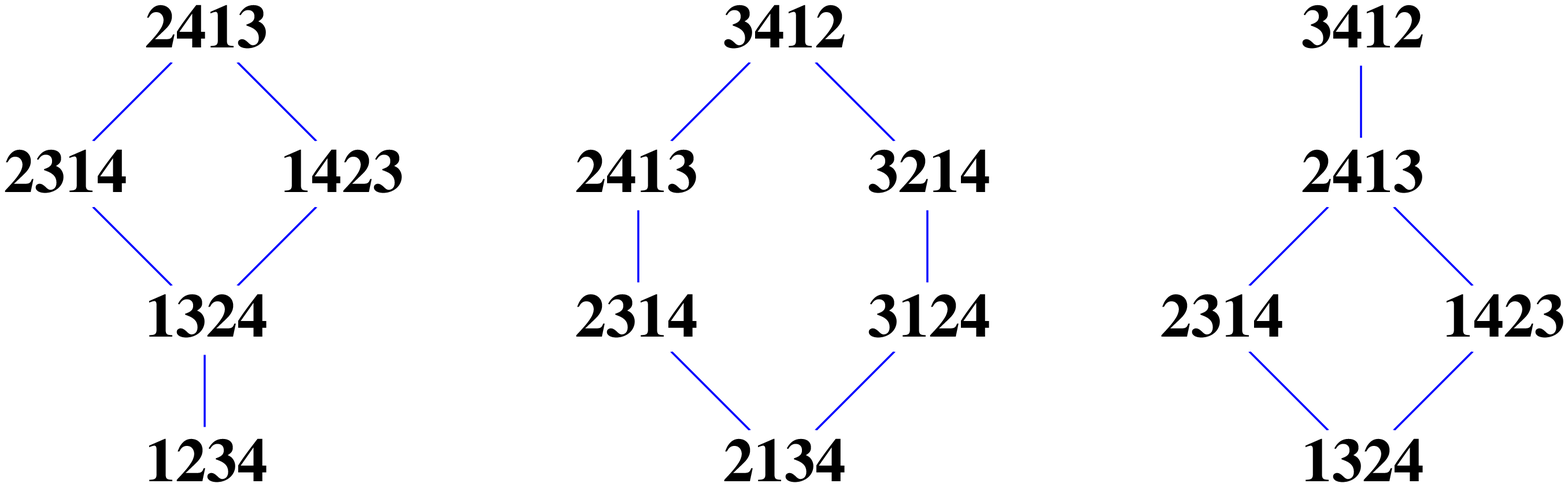}$$

\section{Substitutions}
The identities of \S\ref{sec:Schur_identities}
require a more general study of the behaviour of Schubert polynomials
under certain specializations of the variables.
This leads to a number of new formulas and identities.

For $w\in{\mathcal S}_{n+1}$ and $1\leq p\leq n+1$, let 
$w/_p\in {\mathcal S}_n$ be defined by deleting the $p$th row and $w(p)$th
column from the permutation matrix of  $w$.
If $y\in{\mathcal  S}_n$ and $1\leq q\leq n+1$, then 
$\varepsilon_{p,q}(y)\in {\mathcal S}_{n+1}$ is the permutation such that
$\varepsilon_{p,q}(y)/_p = y$ and $\varepsilon_{p,q}(y)(p)=q$.
The index of summation in a particular case of the Pieri-type
formula~\cite{Lascoux_Schutzenberger_polynomes_schubert,%
sottile_pieri_schubert,Winkel_multiplication}, 
$$
{\mathfrak S}_v \cdot (x_1\cdots x_{p-1}) \quad =\quad
\sum_{v \cpp w} {\mathfrak S}_w,
$$
defines the relation $v\cpp w$.
More concretely,  $v\cpp w$ if and only if there is a chain in the 
$(p-1)$-Bruhat order:
$$
v\ \lessdot_{p-1}\ (\alpha_1,\beta_1) v\ \lessdot_{p-1}\  \cdots\ 
\lessdot_{p-1}\ 
 (\alpha_{p-1},\beta_{p-1})\cdots (\alpha_1,\beta_1) v
\ =\ w
$$
with $\beta_1>\beta_2>\cdots>\beta_{p-1}$.

Define $\Psi_p:{\mathbb Z}[x_1,x_2,\ldots] \rightarrow
{\mathbb Z}[x_1,x_2,\ldots]$ by 
$$
\Psi_p(x_j)\ =\ \left\{\begin{array}{ll} x_j&\mbox{ if } j<p\\
0& \mbox{ if } j=p\\
x_{j-1}&\mbox{ if } j>p\end{array}\right..
$$

\begin{thm}\label{thm:A}
Let $u,w\in {\mathcal S}_\infty$.
Suppose $w(p)=u(p)$ and 
$\ell(w)-\ell(u)=\ell(w/_p)-\ell(u/_p)$, for some positive integer $p$.
Then
\begin{enumerate}
\item[({\em i})]
$\varepsilon_{p,u(p)} :
[u/_p,w/_p]  \stackrel{\sim}{\longrightarrow}  [u,w]$.
\item[({\em ii})] For any $v\in {\mathcal S}_\infty$,
$$
	c^w_{u\, v}\quad =\quad
	\sum_{\stackrel{\mbox{\scriptsize $y\in{\mathcal S}_\infty$}}%
	{v \cpp \varepsilon_{p,1}(y)}} c^{w/_p}_{u/_p\: y}.
$$
\item[({\em iii})] 
For any $v\in {\mathcal S}_\infty$,
$$
 \Psi_p({\mathfrak S}_v) =  
\sum_{\stackrel{\mbox{\scriptsize $y\in{\mathcal S}_\infty$}}%
	{v \cpp \varepsilon_{p,1}(y)}} {\mathfrak S}_y.
$$
\end{enumerate}
\end{thm}

The first statement  is
proven  using combinatorial
arguments, while the second and third
are proven by computing certain maps on cohomology. 
Since  $c^w_{u\,v}=c^w_{v\,u}=c^{w_0 u}_{v\:w_0w}$, 
Theorem~\ref{thm:A} ({\em ii}) gives a
recursion for $c^w_{u\, v}$ when one of $wu^{-1}, wv^{-1}$, or $w_0uv^{-1}$
has a fixed point and the condition on lengths is satisfied.

Theorem~\ref{thm:A} ({\em iii}) is both a
generalization and a strengthening of the transition equations
of~\cite{Lascoux_Schutzenberger_schub_LR_rule}. 
We give an example:
Let $v=413652$.   
Consider the part of the labeled Hasse diagram in the 2-Bruhat order
above $v$ with decreasing  edge labels.
Then the leaves are those $w$ with $v\stackrel{c_3}{\longrightarrow}w$.
$$
\epsfxsize=2.2in \epsfbox{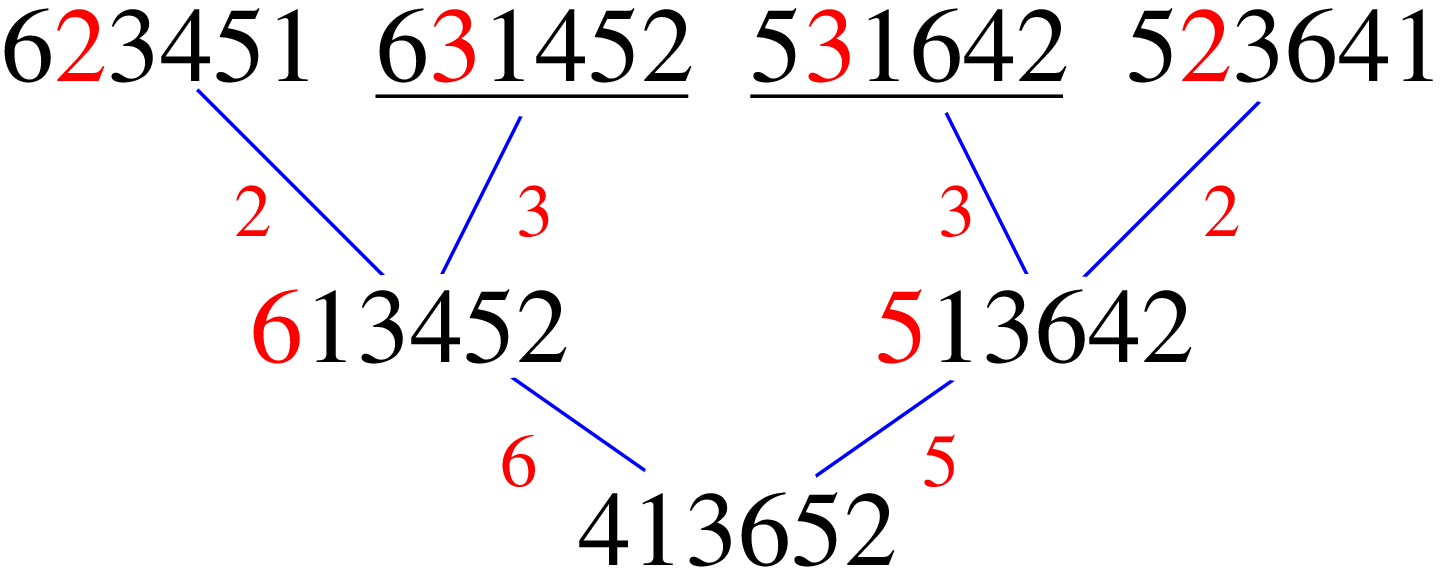}
$$
Of these, only the two underlined permutations are of the form
$\varepsilon_{3,1}(y)$: 
$$
631452\ =\ \varepsilon_{3,1}(52341)\quad\mbox{and}\quad
531642\ =\ \varepsilon_{3,1}(42531).
$$
Thus Theorem~\ref{thm:A} ({\em iii}) asserts that 
$\Psi_3({\mathfrak S}_{413652}) = 
{\mathfrak S}_{52341} + {\mathfrak S}_{42531}.$
Indeed, 
\begin{eqnarray*}
{\mathfrak S}_{413652} &=&
x_1^4x_2x_4x_5 + x_1^3x_2^2x_4x_5 + x_1^3x_2x_4^2x_5 +\\
&\ &
x_1^4x_2x_3x_4 + x_1^4x_2x_3x_5 + x_1^4x_3x_4x_5 + 
x_1^3x_2^2x_3x_4 + x_1^3x_2^2x_3x_5 + x_1^3x_2x_3^2x_4 +\\
&\ & 
x_1^3x_2x_3^2x_5 + x_1^3x_2x_3x_4^2 + 
x_1^3x_3^2x_4x_5 + x_1^3x_3x_4^2x_5 + 
2\cdot x_1^3x_2x_3x_4x_5,
\end{eqnarray*}
so
$$
\Psi_3({\mathfrak S}_{413652}) =
x_1^4x_2x_3x_4 + x_1^3x_2^2x_3x_4 + x_1^3x_2x_3^2x_4.
$$
Since
$$ %\begin{eqnarray*}
{\mathfrak S}_{52341}\ =\ x_1^4x_2x_3x_4\quad\mbox{ and}\quad
{\mathfrak S}_{42531}\ =\  x_1^3x_2^2x_3x_4 + x_1^3x_2x_3^2x_4,
$$%\end{eqnarray*}
we see that 
$$
\Psi_3({\mathfrak S}_{413652}) =
{\mathfrak S}_{52341} + {\mathfrak S}_{42531}.
$$

We also compute the effect of other substitutions of the variables in
terms of the Schubert basis: 
For $P\subset {\mathbb N}$,  set $P^c:= {\mathbb N}-P$ and list the
elements of $P$ and $P^c$ in order:
$$
P\ :\ p_1<p_2<\cdots\qquad \qquad
P^c\ :\ p^c_1<p^c_2<\cdots
$$
Define $\Psi_P:{\mathbb Z}[x_1,x_2,\ldots]\rightarrow
{\mathbb Z}[y_1,y_2,\ldots, z_1,z_2,\ldots]$ by:
$$
\Psi_P(x_{p_j}) \ =\ y_j\qquad\mbox{and}\qquad
\Psi_P(x_{p^c_j}) \ =\ z_j.
$$

\begin{thm}\label{thm:substitution}
Let $P\subset {\mathbb N}$ be as above.
Then there exists an (explicitly described) infinite set 
$\Pi_P\subset {\mathcal S}_\infty$ such that for any $w\in {\mathcal S}_n$ 
and  $\pi\in I_P - {\mathcal S}_{3n}$,
$$
\Psi_P({\mathfrak S}_w)\ =\ 
\sum_{u,\, v} 
c^{(u\times v)\cdot \pi}_{\pi\; w}\;
{\mathfrak S}_u(y)\;{\mathfrak S}_v(z).
$$
\end{thm}

This generalizes~\cite[1.5]{Lascoux_Schutzenberger_structure_de_Hopf}
(See also~\cite[4.19]{Macdonald_schubert}), where it is shown that
the coefficients are nonnegative when $P=[n]$.
Theorem~\ref{thm:substitution} gives infinitely many identities of
the form 
$c^{(u\times v)\cdot \pi}_{\pi\; w} =
c^{(u\times v)\cdot \sigma}_{\sigma\; w}$ for $\pi,\sigma\in \Pi_P$.
Moreover, for these $u,v$, and $\pi$, we have 
$[\pi,\; (u\times v)\cdot \pi]\simeq [e,u]\times[e,\,v]$, 
which is suggestive of a chain-theoretic basis for these identities.

\section{Outline of geometric proofs}

Many of these results are proven with arguments from geometry.
Our main technique is as follows:
If $u,w\in {\mathcal S}_n$, then 
$$
c^w_{u\,v(\lambda,k)}\ =\ \#\left(
X_{w_0w}\bigcap X'_u \bigcap X''_{v(\lambda,k)}\right),
$$
where $X_{w_0w}$, $X'_u$, and $X''_{v(\lambda,k)}$ are Schubert varieties in
general position in the manifold ${\mathbb F}\ell_n$ of complete flags in
${\mathbb C}^n$.
We reduce this to a computation in $\mbox{\em Grass}(k,n)$, the Grassmann
manifold of $k$-planes in ${\mathbb C}^n$.

Let $\pi_k:{\mathbb F}\ell_n \rightarrow\mbox{\em Grass}(k,n)$ be the 
projection
that sends a complete flag to its $k$-dimensional subspace. 
Since $X''_{v(\lambda,k)}=\pi_k^{-1}(\Omega''_\lambda)$, where 
$\Omega''_\lambda$ is a Schubert subvariety of $\mbox{\em Grass}(k,n)$,
we have
$$
c^w_{u\,v(\lambda,k)}\ =\ \# \pi_k\left(
X_{w_0w}\bigcap X'_u\right) \bigcap \Omega''_\lambda.
$$
Thus it suffices to study 
$\pi_k\left(X_{w_0w}\bigcap X'_u\right)\subset \mbox{\em Grass}(k,n)$,
equivalently, its fundamental cycle in homology, as 
$$
\left[\pi_k\left(X_{w_0w}\bigcap X'_u\right)\right] \ =\ 
\sum_\lambda c^w_{u\,v(\lambda,k)}S_{\lambda^c},
$$
where $S_{\lambda^c}$ is the homology class dual to the fundamental cycle of
$\Omega''_\lambda$.

To prove Theorems 3.1 ({\em ii}) and Theorem 3.5, we first use Theorem 4.1
({\em ii}) to 
reduce to the case of $k=l$ and $wu^{-1}=zx^{-1}$.
Then we explicitly compute a dense subset of 
$X_{w_0w}\bigcap X'_u$ and its image, $Y_{w,u}$, in $\mbox{\em Grass}(k,n)$.
This analysis shows that, up to the action of the general linear group,
$Y_{w,u}$ depends only upon $wu^{-1}$, up to conjugation by $(12\ldots n)$,
whenever $wu^{-1}\in {\mathcal S}_n$.

For Theorem 4.2, we study maps
$$
\Psi_P\ :\ {\mathbb F}\ell_n\times {\mathbb F}\ell_m \ 
\longrightarrow {\mathbb F}\ell_{n+m}
$$
where $\Psi_P$ `shuffles' pairs of flags together to obtain a longer flag,
according to a set $P: 1\leq p_1<\cdots<p_n\leq m$.
We show that $\Psi_P(X_u\times X_v)$ is an intersection of two Schubert
varieties, which enables the computation of the map $(\Psi_P)_*$ on homology.
Then Poincar\'e duality determines the map $(\Psi_P)^*$ on the Schubert
basis of cohomology. 
By construction, $(\Psi_P)^*$ acts by the substitution $\Psi_P$ of \S 4.
In the case $m=1$, these computations become more precise, and we obtain
Theorem 4.1 ({\em i}) and ({\em ii}).

Finally, for Theorem 3.4, 
suppose $\zeta\cdot\eta$ is disjoint and $u\leq_{k+l}(\zeta\cdot\eta)u$
with $u\in {\mathcal S}_{n+m}$.
Then, set $P=u^{-1}{\rm supp}(\zeta)$ and consider the commutative diagram:
$$
\setlength{\unitlength}{2.2pt}
\begin{picture}(107,34)
\put(0,3){$\Gr_k{\mathbb C}^n\times\Gr_l{\mathbb C}^m$}
\put(75,3){$\Gr_{k+l}{\mathbb C}^{n+m}$}
\put(13.8,26){${\mathbb F}\ell_n\times{\mathbb F}\ell_m$}
\put(82,26){${\mathbb F}\ell_{n+m}$}
\put(58,7){$\varphi_{k,l}$}     \put(60,30){$\psi_P$}
\put(6.5,16){$\pi_k\times\pi_l$}  \put(91,16){$\pi_{k+l}$}
\put(53,4.5){\vector(1,0){20}}
\put(40,27.5){\vector(1,0){40}}
\put(25.5,23){\vector(0,-1){14}}
\put(89,23){\vector(0,-1){14}}
\end{picture}
$$
Here, $\varphi_{k,l}$ maps a pair 
$(H,K)\in \Gr_k{\mathbb C}^n\times\Gr_l{\mathbb C}^m$ to the sum 
$H\oplus K$ in $\Gr_{k+l}{\mathbb C}^{n+m}$.
We show there exists $x\in{\mathcal S}_n$, 
$y\in {\mathcal S}_m$ and $\zeta',\eta'$
shape-equivalent to $\zeta$ and $\eta$ such that
$x \leq_k \zeta' x$, $y \leq_l \eta' y$, and
$$
\Psi_P\left(X_{w_0\zeta'x}\bigcap X'_x\right)\times
\left(X_{w_0\eta'y}\bigcap X'_y\right)\ =\ 
X_{w_0(\zeta\cdot\eta)u}\bigcap X'_u.
$$
Thus to compute $\pi_{k+l}(X_{w_0(\zeta\cdot\eta)u}\bigcap X'_u)$,
or rather its homology class, it suffices to compute the map 
$(\varphi_{k,l})_*$ on homology, which is
$$
(\varphi_{k,l})_* S_\lambda\ =\ \sum_{\mu,\nu}c^\lambda_{\mu\,\nu}
S_\mu\otimes S_\nu.
$$
For more details on these proofs and other aspects of this note, 
see~\cite{bergeron_sottile_symmetry}.

\end{document}